\documentclass[MSNbibl,number,citesort,dvips]{arxbj}

\aid{0}
\volume{19}
\issue{1}
\pubyear{2013}
\firstpage{137}
\lastpage{153}
\doi{10.3150/11-BEJ396}

\makeatletter
\newcommand{\eqref}[1]{(\ref{#1})}
\newcommand{\C}{\mathcal{C}}
\newcommand{\F}{\mathcal{F}}
\newcommand{\G}{\mathcal{G}}
\newcommand{\R}{\mathbb{R}}
\renewcommand{\u}{\bar u}
\newcommand{\eqdef}{\stackrel{\mathrm{def}}=}
\newcommand{\pto}{\stackrel{p}\to}
\newtheorem{thmm}{Theorem}[section]
\newtheorem{cor}[thmm]{Corollary}
\newtheorem{lem}[thmm]{Lemma}
\makeatother

\begin{document}
\begin{frontmatter}

\title{Improving Brownian approximations for boundary crossing problems}
\runtitle{Improving Brownian approximations}
\begin{aug}
\author{\fnms{Robert} \snm{Keener}\corref{}\ead[label=e1]{keener@umich.edu}}
\runauthor{R. Keener} 
\address{Department of Statistics,
University of Michigan,
Ann Arbor, MI 48103, USA.\\
\printead{e1}}
\end{aug}

\received{\smonth{7} \syear{2010}}
\revised{\smonth{7} \syear{2011}}

%
\begin{abstract}
Donsker's theorem shows that random walks behave like Brownian
motion in an asymptotic sense. This result can be used to
approximate expectations associated with the time and
location of a random walk when it first crosses a nonlinear
boundary. In this paper, correction terms are derived to improve
the accuracy of these approximations.
\end{abstract}


\end{frontmatter}
%

\section{Introduction and main results}

Let $X,X_1,X_2,\ldots$ be i.i.d.
with mean zero and unit variance; take $S_k=X_1+\cdots+X_k$, $k\ge1$,
with $S_0=0$; and let $W(t)$, $t\ge0$, be standard Brownian
motion. By Donsker's theorem, if $W_n$ is continuous and piecewise
linear with
\[
W_n(k/n)=S_k/\sqrt n, \qquad k=0,1,\ldots,
\]
then $W_n\Rightarrow W$ in $\C[0,\infty)$ as $n\to\infty$.
Let $b$ be a smooth function on $[0,\infty)$ with $b(0)>0$,
such that
\[
\tau_0=\inf\{t\ge0\dvt W(t)\ge b(t)\}
\]
is finite almost surely, and define
\[
\tau_n=\inf\{k/n\ge0\dvt W_n(k/n)\ge b(k/n)\}.
\]
Defining boundary levels
%
\begin{equation}b_k=b_{k,n}=\sqrt nb(k/n)
\label{blevels},
\end{equation}
this stopping time can be written as
\[
\tau_n=\inf\{k\ge0\dvt S_k\ge b_k\}/n.
\]
As the form suggests, $\tau_n\Rightarrow\tau_0$ as $n\to\infty$.
This can be established by introducing $\tilde\tau_n
=\inf\{t\ge0\dvt W(t)\ge b(t)\}$, arguing that
$\tau_n-\tilde\tau_n\pto0$, and using the continuous mapping
theorem, Theorem~{5.1} of Billingsley~\cite{bi68}, to show that $\tilde\tau_n
\Rightarrow\tau_0$. Note that the Brownian path $W(\cdot)$ will
be a continuity point for the transformation $W(\cdot)\leadsto\tau_0$
whenever $\tau_0=\inf\{t\ge0\dvt W(t)>b(t)\}$, and this holds
with probability one by the strong Markov property. We also have
$W_n(\tau_n)-b(\tau_n)\pto0$, and so
\[
(\tau_n,W_n(\tau_n))
\Rightarrow(\tau_0,W(\tau_0))
\]
as $n\to\infty$. Thus if $f$ is a bounded continuous function,
%
\begin{equation}\label{thelimit}
Ef(\tau_n,W_n(\tau_n))\to Ef(\tau_0,W(\tau_0)).
\end{equation}
For large $n$, the limit here provides a natural approximation for the
expectation on the left-hand side. The main result of this paper provides
correction terms of order $1/\sqrt n$, improving this approximation.
The excess over the boundary,
\[
R_n=S_{n\tau_n}-\sqrt nb(\tau_n)
=\sqrt n[W_n(\tau_n)-b(\tau_n)],
\]
plays an important role in this analysis. The excess over the boundary
also plays a central role in nonlinear renewal theory, where the law
of large numbers drives the leading order approximation. See
Woodroofe~\cite{wo82} or Siegmund~\cite{si85} for a discussion and applications to
sequential analysis. With the Brownian motion scaling considered in
this paper, results on improved approximations and the excess over the
boundary are given by Siegmund~\cite{si79}, Siegmund and Yuh~\cite{sy82}, Yuh~\cite{yu82} and
Hogan~\cite{ho84a,ho84b,ho84c,ho86}. Siegmund~\cite{si85} suggests various
applications of this theory to sequential analysis; Broadie
\textit{et~al.}~\cite{bgk97},
Broadie
\textit{et~al.}~\cite{bgk99} and Kou~\cite{ko03} use it to study options pricing; and
Glasserman and Liu~\cite{gl97} consider its use for inventory control. With the
exception of Hogan~\cite{ho84b,ho84c}, stopping boundaries in these
papers are linear.

To appreciate the role of the excess $R_n$ in improving
\eqref{thelimit}, note that if
$f(t,x)=h(t)[x-b(t)]$, then $Ef(\tau_n,W_n(\tau_n))
=Eh(\tau_n)R_n/\sqrt n$. Hogan~\cite{ho84b} derives the
limiting joint distribution for $R_n$ and $\tau_n$; they are asymptotically
independent, and the limiting distribution for $R_n$ has mean
%
\begin{equation}
\rho=\frac{ES_{T_0}^2}{2ES_{T_0}},
\label{rhodef}
\end{equation}
where $T_0$ is the ladder time
\[
T_0=\inf\{k>0\dvt S_k\ge0\}.
\]
Hogan's argument is quite delicate. It is based on conditioning on a
stopping time with a boundary just slightly less than the boundary for
$\tau_n$. By contrast, the approach pursued here is more global and
analytic in character, but relies on smoothness of $f$ and $b$ to a
greater extent. Formulas to calculate $\rho$ numerically are given by
Siegmund~\cite{si85} and Keener~\cite{ke94}.

An important special case of \eqref{thelimit} would be first passage
probabilities, $P(\tau_n\le t)$. The regularity conditions here
require differentiable $f$, so this case is formally excluded
(although our result would suggest an approximation). Refined
approximations for these probabilities are also suggested by
Hogan~\cite{ho84b}, but his derivation is heuristic and assumes $EX^3=0$.

The limit in \eqref{thelimit} can be found by solving the heat
equation. To
describe its relevance, let $Y=Y(t,x)$ be a process starting at time
$t$ and position $x$ given by
\[
Y_s=Y_s(t,x)=x+W(s-t),\qquad  s\ge t;
\]
let $\tau=\tau(t,x)$ be stopping times given by
\[
\tau=\tau(t,x)=\inf\{s\ge t\dvt Y_s\ge b(s)\};
\]
and define
\[
u(t,x)=Ef(\tau,Y_\tau), \qquad t\ge0, x\le b(t).
\]
Noting that $\tau_0=\tau(0,0)$ and $W(\tau_0)=Y_{\tau(0,0)}$, the
limit $Ef(\tau_0,W(\tau_0))$ in \eqref{thelimit} is $u(0,0)$. By
the Feynman--Kac formula (Kac~\cite{ka51}), $u$ satisfies the heat
equation
\[
u_t+\frac12u_{xx}=0
\]
in the region $\{(t,x)\dvt t\ge0,x<b(t)\}$, with boundary condition
$u(t,b(t))=f(t,b(t))$. Furthermore, $u$ is the
unique solution in a suitable class of functions; see Krylov~\cite{kr91} or
Bass~\cite{ba98}. In practice, $u(0,0)$ can be computed by numerical
solution of the heat equation. In the sequel, continuity and
differentiability of $u$ will play an important role.

Boundary effects associated with the excess $R_n$ only arise (to order
$\mathrm{o}(1/\sqrt n)$) when $f_x$ and $u_x$ disagree along the boundary. Let
$\Delta(t)$ denote the difference
\[
\Delta(t)=f_x(t,b(t))-u_x(t,b(t)-),\qquad  t>0,
\]
and decompose $f$ as the sum $f_0+f_1$ with
\[
f_0(t,x)=f(t,x)-\Delta(t)\bigl(x-b(t)\bigr)
\]
and
\[
f_1(t,x)=\Delta(t)\bigl(x-b(t)\bigr).
\]
Since $u$ and $u_x$ agree with $f_0$ and $\partial f_0/\partial x$
along the boundary, it seems appropriate to view $u(0,0)$ as an
approximation for $Ef_0(\tau_n,W_n(\tau_n))$. It is then
natural and convenient to extend $u$ above the boundary, defining
\[
\u(t,x)=
\cases{ u(t,x),&\quad $x\le b(t);$\vspace*{2pt}\cr
f_0(t,x),&\quad $x>b(t).$}
\]
With this convention, $\u$ and $\u_x$ are both continuous at the boundary.
Note also that
\[
Ef_1(\tau_n,W_n(\tau_n))=\frac{1}{\sqrt n}ER_n\Delta(\tau_n).
\]

\begin{thmm}\label{thm:main}
Assume:
\begin{enumerate}
\item The distribution of $X$ is strongly non-lattice (or satisfies Cramer's
condition C),
\[
\limsup_{|t|\to\infty}|E\mathrm{e}^{\mathrm{i}tX}|<1,
\]
and $EX=0$, $EX^2=1$ and $EX^4<\infty$.
\item The stopping times $\tau_n$, $n\ge1$, are uniformly integrable.
\item The boundary function $b$ has a bounded first derivative and $b(0)>0$.
\item The function $f$ and its first and second order partial derivatives
are bounded and continuous.
\item The functions $u$, $u_x$, $u_{xx}$, $u_{xxx}$, $u_{xxxx}$,
$u_t$ and $u_{tt}$ are bounded and continuous.
\end{enumerate}
Then
\begin{eqnarray*}
Ef(\tau_n,W_n(\tau_n))
&=&Ef(\tau_0,W(\tau_0))
+\frac{EX^3}{6\sqrt n}
E\int_0^{\tau_0}u_{xxx}(t,W(t)) \,\mathrm{d}t\\
&&{} +\frac{\rho}{\sqrt n}E\Delta(\tau_0)+\mathrm{o}\bigl(1/\sqrt n\bigr)
\end{eqnarray*}
as $n\to\infty$.
\end{thmm}

The second assumption will hold if $b(t)+\epsilon t\to-\infty$ for
some $\epsilon>0$. If $b$ and $f$ are sufficiently smooth, then the
final assumption follows from standard H\"older estimates for
solutions of parabolic differential equations; see, for instance,
Problem~4.5 of Lieberman~\cite{li96}.

The heat equation for $u$ can be derived, at least informally, by
conditioning a short time interval into the future. There is an analogous
equation in discrete time. Define
\[
\tau_n(t,x)=\inf\bigl\{t+k/n\dvt x+S_k/\sqrt n\ge b(t+k/n),k=0,1,\ldots
\bigr\}
\]
and
\[
u_n(t,x)=Ef_0\bigl(\tau_n(t,x),x+S_{n\tau_n(t,x)}/\sqrt n\bigr).
\]
Conditioning on $X_1$,
%
\begin{equation}\label{convolution}
u_n(t,x)=
\cases{ f_0(t,x),&\quad $x\ge b(t);$\vspace*{2pt}\cr
Eu_n\bigl(t+1/n,x+X/\sqrt n\bigr),&\quad $x<b(t).$}
\end{equation}
Unfortunately, with integration against the distribution of $X$, this
convolution-type equation is usually less tractable numerically than
the heat equation.

Theorem~\ref{thm:main} evolved from my attempts to improve $\u$ as an
approximation for $u_n$ by imitating the matched asymptotic expansions
used to study boundary effects in partial\vadjust{\goodbreak} differential equations. The
method might also be viewed as martingale approximation, with bounds
for potential or renewal measures playing a central role in the proofs.

To study the error of $u(0,0)=\u(0,0)$ as an approximation for
$Ef_0(\tau_n,W_n(\tau_n))$, define functions
%
\begin{equation}\label{edef}
e_n(t,x)=
\cases{
E\u\bigl(t+1/n,x+X/\sqrt n\bigr)-\u(t,x),&\quad $x<b(t);$\vspace*{2pt}\cr
0,&\quad $x\ge b(t).$}\
\end{equation}
Writing
%
\begin{eqnarray}\label{telescope}
Ef_0(\tau_n,W_n(\tau_n))
&=&E\u((\tau_n,W_n(\tau_n))\nonumber\\
&=&\u(0,0)+E\sum_{k=0}^{n\tau_n-1}\bigl[\u\bigl((k+1)/n,S_{k+1}/\sqrt
n\bigr)
-\u\bigl(k/n,S_k/\sqrt n\bigr)\bigr]\\
&=&\u(0,0)+E\sum_{k=0}^{n\tau_n-1}e_n\bigl(k/n,S_k/\sqrt
n\bigr),\nonumber
\end{eqnarray}
a correction term for the approximation $\u(0,0)$ will be sought by
approximating the expected sum in this equation. Details for
this calculation are given in Section~\ref{sec:approxf0}.
The approximation for $Ef_1(\tau_n,W_n(\tau_n))$
is derived in Section~\ref{sec:approxf1}.

\section{\texorpdfstring{An approximation for $Ef_0(\tau_n,W_n(\tau_n))$}
{An approximation for Ef 0(tau n,W n(tau n))}}
\label{sec:approxf0}

\begin{lem}\label{lem:approxe}
Under the assumptions of Theorem~\ref{thm:main},
%
\begin{eqnarray}\label{firstline}
E\u\bigl(t,x+X/\sqrt n\bigr)
&=&\u(t,x)+\frac{\u_{xx}(t,x)}{2n}+\frac{EX^3\u_{xxx}(t,x)}{6n\sqrt n}
\nonumber
\\[-8pt]
\\[-8pt]
\nonumber
&&{} +\mathrm{O}(1/n^2)+\mathrm{O}\biggl(\frac{1/n}{1+n[b(t)-x]^2}\biggr)
\end{eqnarray}
as $n\to\infty$, uniformly for $t\ge0$, $x<b(t)$. From this,
\[
e_n(t,x)=\frac{EX^3\u_{xxx}(t,x)}{6n\sqrt n}+\mathrm{O}(1/n^2)
+\mathrm{O}\biggl(\frac{1/n}{1+n[b(t)-x]^2}\biggr)
\]
as $n\to\infty$, uniformly for $t\ge0$, $x<b(t)$.
\end{lem}

\begin{pf}
By Taylor expansion of $u$, on $\{x+X/\sqrt n\le b(t)\}$ we have
%
\begin{eqnarray}\label{below}
\u\bigl(t,x+X/\sqrt n\bigr)
&=&\u(t,x)+\frac{X\u_x(t,x)}{\sqrt n}+\frac{X^2\u_{xx}(t,x)}{2n}
\nonumber
\\[-8pt]
\\[-8pt]
\nonumber
&&{} +\frac{X^3\u_{xxx}(t,x)}{6n\sqrt n}+\mathrm{O}(X^4/n^2).
\end{eqnarray}
Lagrange's formula for the remainder will involve $\u_{xxxx}$ at an
intermediate value $x^*$ between $x$ and $x+X/\sqrt n$, and from this
it is clear that this equation holds uniformly for $t\ge0$, $x<b(t)$.
Since $\u_{xx}(t,x)$ exists unless $x=b(t)$ and is bounded, on
$\{x+X/\sqrt n>b(t)\}$,
%
\begin{equation}\label{above}
\u\bigl(t,x+X/\sqrt n\bigr)
=\u(t,x)+\frac{X\u_x(t,x)}{\sqrt n}+\mathrm{O}(X^2/n),
\end{equation}
as $n\to\infty$. Again, this will hold uniformly for $t\ge0$, $x<b(t)$.
Noting that
\[
\frac{|X|^3}{n\sqrt n}\le\frac{X^2}{n}+\frac{X^4}{n^2},
\]
we can combine \eqref{below} and \eqref{above} to obtain
\begin{eqnarray*}
u\bigl(t,x+X/\sqrt n\bigr)
&=&u(t,x)+\frac{Xu_x(t,x)}{\sqrt n}+\frac{X^2u_{xx}(t,x)}{2n}
+\frac{X^3u_{xxx}(t,x)}{6n\sqrt n}\\
&&{} +\mathrm{O}(X^4/n^2)+\mathrm{O}(X^2/n)I\bigl\{X>\sqrt n\bigl(b(t)-x\bigr)\bigr\}.
\end{eqnarray*}
The first assertion \eqref{firstline} follows by integrating against the
distribution of $X$, noting that
\begin{eqnarray*}
E\bigl[X^2;X>\sqrt n\bigl(b(t)-x\bigr)\bigr]
&\le&\min\biggl\{\frac{EX^4}{n[b(t)-x]^2},EX^2\biggr\}\\
&\le&\frac{1+EX^4}{1+n[b(t)-x]^2}.
\end{eqnarray*}
Here and in the sequel, $E[Y;B]\eqdef E(Y1_B)$.

If $x<b(t)$ and $x<b(t+1/n)$, then, by \eqref{edef} and \eqref{firstline},
\begin{eqnarray*}
e_n(t,x)&=&\u(t+1/n,x)-u(t,x)+\frac{\u_{xx}(t+1/n,x)}{2n}
+\frac{EX^3\u_{xxx}(t+1/n,x)}{6n\sqrt n}\\
& &{}+\mathrm{O}(1/n^2)+\mathrm{O}\biggl(\frac{1/n}{1+n[b(t+1/n)-x]^2}\biggr).
\end{eqnarray*}
In this case, the second assertion follows by the Taylor expansion
\begin{eqnarray*}
\u(t+1/n,x)
&=&\u(t,x)+\frac1n\u_t(t,x)+\mathrm{O}(1/n^2)\\
&=&\u(t,x)-\frac{1}{2n}\u_{xx}(t,x)+\mathrm{O}(1/n^2),
\end{eqnarray*}
and because
\[
\frac{1+n[b(t)-x]^2}{1+n[b(t+1/n)-x]^2}
\]
is uniformly bounded as $b'$ is bounded. If, instead,
$x\ge b(t+1/n)$, but $x<b(t)$, then $n[b(t)-x]^2\to0$
and $n[b(t+1/n)-x]^2\to0$, and the asymptotic bound holds
because $u(t+1/n,x)-u(t,x)=\mathrm{O}(1/n)$.
\end{pf}

%

Define
\[
N_d=N_d(n)=\#\{k<n\tau_n\dvt S_k>b_k-d\},
\]
the number of times the walk is within distance $d$ of the boundary
before stopping. The following result is essentially due to Hogan~\cite{ho84b}.
It slightly improves a bound given in the proof for Lemma~1.1 in his paper.

\begin{lem}\label{lem:edge}
With the assumptions of Theorem~\ref{thm:main}, there exists a finite
constant $K\ge0$ such that
\[
EN_d=K(1+d^2),
\]
for all $n\ge1$ and $d>0$.
Also, if
\[
M_B(\alpha)=\#\{k\le\alpha n\tau_n\dvt b_k-S_k\in B\},
\]
then there exists a finite constant $K>0$ such that
\[
EM_B(\alpha)\le KP\bigl(M_B(\alpha)\ge1\bigr)\bigl(1+(\sup B)^2\bigr),
\]
for all $n\ge1$, all $\alpha>0$ and all $B\subset\R$.
\end{lem}

\begin{pf}
Without loss of generality, let $d$ be a positive integer.
By the central limit theorem,
%
\begin{equation}\label{geobound}
P\{S_{n^2}>(1+\Vert b'\Vert_\infty)n\}\ge\gamma>0,
\end{equation}
for all $n$ sufficiently large, say $n\ge n_0$. Since
the $\tau_n$ are uniformly integrable (by the second assumption of Theorem~\ref{thm:main}) and $N_d\le
n\tau_n$,
we can assume that $n_0d^2\le n$. Define
\[
N_{m,d}=\#\{k\le m\dvt k<n\tau_n,S_k>b_k-d\},
\]
and let
\[
\nu_{j,d}=\inf\{m\dvt N_{m,d}=j\},
\]
so the $j$th time the walk is within $d$ of the boundary happens on
step $\nu_{j,d}$. Note that $N_d\ge j+n_0d^2$ implies the walk is below
the boundary at time $\nu_{j,d}+n_0d^2$, that is,
\[
S_{\nu_{j,d}+n_0d^2}<b_{\nu_{j,d}+n_0d^2},
\]
which, in turn, implies
\[
S_{\nu_{j,d}+n_0d^2}-S_{\nu_{j,d}}<b_{\nu_{j,d}+n_0d^2}-b_{\nu_j}+d
\le d+n_0d^2\frac{\Vert b'\Vert_\infty}{\sqrt n}
\le\sqrt{n_0}d(1+\Vert b'\Vert_\infty).
\]
But $S_{\nu_{j,d}+n_0d^2}-S_{\nu_{j,d}}$ is independent of $\{N_d\ge k\}$.
So using this bound and \eqref{geobound},
\[
P(N_d\ge j+n_0d^2)\le P(N_d\ge j)(1-\gamma).
\]
Iterating this,
\begin{eqnarray*}
&& P(N_d\ge1+jn_0d^2)\\
&&\quad =P\bigl(N_d\ge1+(j-1)n_0d^2+n_0d^2\bigr)\\
&&\quad\le P\bigl(N_d\ge1+(j-1)n_0d^2\bigr)(1-\gamma)
\le\cdots
\le P(N_d\ge1)(1-\gamma)^j, \qquad j=0,1,\ldots.
\end{eqnarray*}
Hence
\begin{eqnarray*}
EN_d
&=&\int_0^\infty P(N_d\ge x) \,\mathrm{d}x\\
&\le& P(N_d\ge1)\biggl[1+\int_1^\infty
(1-\gamma)^{\lfloor(x-1)/(n_0d^2)\rfloor} \,\mathrm{d}x\biggr]\\
&=&P(N_d\ge1)\biggl[1+\frac{n_0d^2}{\gamma}\biggr].
\end{eqnarray*}
The proof of the bound for $EM_B(\alpha)$ is the same.
\end{pf}

\begin{cor}\label{cor:sumbounds}
Let $c_k=c_{k,n}$, $k\ge0$, $n\ge1$ be constants. Define
\[
\Lambda=\sup_{k,n}(b_k-c_k),
\]
and let $g$ be a non-negative function on $(-\infty,\Lambda]$.
If $\Lambda<\infty$, $\Vert g\Vert_\infty<\infty$, and
$g(x)\to0$ as $x\to-\infty$,
\[
\frac1nE\sum_{k=0}^{n\tau_n-1}f(S_k-c_k)\to0,
\]
as $n\to\infty$. If, in addition, $g$ is non-decreasing,
\[
E\sum_{k=0}^{n\tau_n-1}g(S_k-c_k)
\le K\biggl[g(\Lambda)+2\int_{-\infty}^0|x|g(x+\Lambda) \,\mathrm{d}x\biggr],
\]
where $K$ is the constant in Lemma~\ref{lem:edge}.
\end{cor}

When this corollary is used later, $c_k$ will be either $b_k$ or
$b_{k+1}$. When $c_k=b_k$, $\Lambda$ is zero, and when $c_k=b_{k+1}$,
$\Lambda\le\Vert b'\Vert_\infty$.

\begin{pf*}{Proof of Corollary \protect\ref{cor:sumbounds}}
For the first assertion, for any $d>0$,
\[
g(S_k-c_k)\le\Vert g\Vert_\infty I\{S_k>b_k-d\}
+\sup_{x\in(-\infty,\Lambda-d]}g(x).
\]
Summing over $k$ and bounding the expectation using Lemma~\ref{lem:edge},
\[
\frac1nE\sum_{k=0}^{n\tau_n-1}g(S_k-c_k)\le
\frac1n\Vert g\Vert_\infty K(1+d^2)
+\sup_{x\in(-\infty,\Lambda-d]}g(x)E\tau_n,
\]
and the result follows because $d$ can be arbitrarily large.

In the second assertion, we can assume without loss of generality
that $g$ is right continuous and write
\[
g(y)=\int I\{x\le y\} \,\mathrm{d}g(x).
\]
By Fubini's theorem and Lemma~\ref{lem:edge},
\begin{eqnarray*}
E\sum_{k=0}^{n\tau_n-1}g(S_k-c_k)
&\le& E\sum_{k=0}^{n\tau_n-1}g(S_k-b_k+\Lambda)\\
&\le&\int E\sum_{k\ge0}I\{x<S_k-b_k+\Lambda,k<n\tau_n\} \,\mathrm{d}g(x)\\
&=&\int EN_{\Lambda-x} \,\mathrm{d}g(x)\\
&\le& K\int_{-\infty}^\Lambda[1+(\Lambda-x)^2] \,\mathrm{d}g(x)\\
&=&K\biggl[g(\Lambda)+2\int_{-\infty}^0|x|g(x+\Lambda) \,\mathrm{d}x\biggr].
\end{eqnarray*}
\upqed\end{pf*}

The second assertion in Corollary~\ref{cor:sumbounds} is useless when
the integral in the bound diverges, but, in certain cases, it gives
sharper results than the first assertion. The next corollary
considers a specific function of interest later.

\begin{cor}\label{cor:edge}
With the assumptions of Theorem~\ref{thm:main},
\[
E\sum_{k=0}^{n\tau_n-1}\frac{1}{1+[b_k-S_k]^2}
=\mathrm{O}(\log n)
\]
as $n\to\infty$.
\end{cor}

\begin{pf}
If $0\le b_k-S_k\le\sqrt n$,
\[
\frac{1}{1+[b_k-S_k]^2}
=\frac{1}{1+n}+\int\frac{I\{b_k-S_k<x<\sqrt n\}2x}{(1+x^2)^2} \,\mathrm{d}x,
\]
and so
\[
\sum_{k=0}^{n\tau_n-1}\frac{1}{1+[b_k-S_k]^2}
\le\tau_n+\int_0^{\sqrt n}\frac{2xN_x}{(1+x^2)^2} \,\mathrm{d}x.
\]
Using Lemma~\ref{lem:edge},
\[
E\sum_{k=0}^{n\tau_n-1}\frac{1}{1+[b_k-S_k]^2}
=E\tau_n+\mathrm{O}(1)\int_0^{\sqrt n}\frac{2x}{1+x^2} \,\mathrm{d}x=\mathrm{O}(\log n).
\]
\upqed\end{pf}

The final corollary gives uniform integrability for moments of $R_n$.

\begin{cor}\label{cor:unifint}
With the assumptions of Theorem~\ref{thm:main}, if $E|X|^{p+2}<\infty$,
$R^p_n$, $n\ge1$, are uniformly integrable.
\end{cor}

\begin{pf}
Conditioning on $\F_k=\sigma(X_1,\ldots,X_k)$, if $c>0$,
\begin{eqnarray*}
E[R_n^p;R_n\ge c]
&=&\sum_{k\ge0}E[(S_k+X_{k+1}-b_{k+1})^p;k<n\tau_n,
S_k+X_{k+1}-b_{k+1}\ge c]\\
&=&E\sum_{k<n\tau_n}g(S_k-b_{k+1}),
\end{eqnarray*}
where
\[
g(x)=E[(x+X)^p;x+X\ge c].
\]
This function is increasing and right continuous. Taking
$\Lambda=\sup_{k,n}(b_k-b_{k+1})\le\Vert b'\Vert_\infty$,
by Fubini's theorem,
\begin{eqnarray*}
\int_{-\infty}^0|x|g(x+\Lambda) \,\mathrm{d}x
&=&-E\int x(x+\Lambda+X)^pI\{c-X-\Lambda\le x<0\} \,\mathrm{d}x\\
&=&E\biggl[\frac{(X+\Lambda)^{p+2}-c^{p+2}}{(p+1)(p+2)}
+\frac{c^{p+1}(X+\Lambda-c)}{p+1};X+\Lambda\ge c\biggr].
\end{eqnarray*}
This expectation tends to zero as $c\to\infty$ by dominated convergence,
as does $g(\Lambda)$, and uniform integrability follows from the bound
in Corollary~\ref{cor:sumbounds}.
\end{pf}

\begin{thmm}\label{thm:f0approx}
Under the assumptions of Theorem~\ref{thm:main},
\[
Ef_0(\tau_n,W_n(\tau_n))=\u(0,0)+\frac{EX^3}{6\sqrt n}
E\int_0^{\tau_0}\u_{xxx}(t,W(t)) \,\mathrm{d}t+\mathrm{o}\bigl(1/\sqrt n\bigr).
\]
\end{thmm}

\begin{pf}
Because
\begin{eqnarray*}
&\displaystyle\frac1n\sum_{k=0}^{n(T\wedge\tau_n)-1}\u_{xxx}\bigl(k/n,S_k/\sqrt n\bigr)
=\int_0^{T\wedge\tau_n}\u_{xxx}\bigl(\lfloor nt\rfloor/n,
W_n(\lfloor nt\rfloor/n)\bigr) \,\mathrm{d}t,&
\\
&\displaystyle\frac1nE\sum_{k=0}^{n(T\wedge\tau_n)-1}\u_{xxx}\bigl(k/n,S_k/\sqrt
n\bigr)
-E\int_0^{T\wedge\tau_n}\u_{xxx}(t,W_n(t)) \,\mathrm{d}t\to0&
\end{eqnarray*}
by dominated convergence, since $\max_{k<nT}|X_k|/\sqrt n\to0$
almost surely. So, by Donsker's theorem,
\[
E\int_0^{T\wedge\tau_n}\u_{xxx}(t,W_n(t)) \,\mathrm{d}t\to
E\int_0^{T\wedge\tau}\u_{xxx}(t,W(t)) \,\mathrm{d}t.
\]
Then, since $u_{xxx}$ is uniformly bounded and $\tau_n$, $n\ge1$ are
uniformly integrable,
\[
\frac1nE\sum_{k=0}^{n\tau_n-1}\u_{xxx}\bigl(k/n,S_k/\sqrt n\bigr)
\to E\int_0^{\tau_0}\u_{xxx}(t,W(t)) \,\mathrm{d}t.
\]
The theorem now follows from \eqref{telescope} using the formula for
$e_n$ in Lemma~\ref{lem:approxe} and the asymptotic bound in
Corollary~\ref{cor:edge}.
\end{pf}

\section{\texorpdfstring{An approximation for $ER_n\Delta_n(\tau_n)$}
{An approximation for ER n Delta n(tau n)}}
\label{sec:approxf1}

\begin{thmm}\label{thm:R}
Under the assumptions of Theorem~\ref{thm:main},
\[
ER_n\Delta(\tau_n)\to\rho E\Delta(\tau_0),
\]
where $\rho$ is the limiting mean excess defined in \eqref{rhodef}.
\end{thmm}

The proof of this result, like that for Theorem~\ref{thm:f0approx},
is based on a telescoping sum argument, but now the summands involve
functions related to fluctuation theory for random walks. For $x\le0$,
define stopping times
\[
T_x=\inf\{k\ge1\dvt x+S_k\ge0\},
\]
and define
\[
H(x)=
\cases{ x-\rho,&\quad $x\ge0;$\vspace*{2pt}\cr
E[S_{T_x}+x]-\rho;&\quad $x<0.$}
\]
Conditioning on $X_1$, for $x<0$
%
\begin{equation}\label{conveqnH}
H(x)=EH(x+X).
\end{equation}
In particular, on $\{S_k<b_k\}$,
%
\begin{equation}\label{Heqn}
E[H(S_{k+1}-b_{k})\big\vert\F_k]=H(S_k-b_k).
\end{equation}
Now
\[
R_n\Delta(\tau_n)=\rho\Delta(\tau_n)
+\Delta(\tau_n)H\bigl[S_{n\tau_n}-\sqrt nb(\tau_n)\bigr],
\]
and by a telescoping sum argument,
\begin{eqnarray*}
&& E\Delta(\tau_n)H\bigl[S_{n\tau_n}-\sqrt nb(\tau_n)\bigr]
-\Delta(0)H\bigl[-\sqrt nb(0)\bigr]\\
&&\quad=E\sum_{k=0}^{n\tau_n-1}\biggl[\Delta\biggl(\frac{k+1}{n}\biggr)
H(S_{k+1}-b_{k+1})-\Delta\biggl(\frac{k}{n}\biggr)
H(S_k-b_k)\biggr]\\
&&\quad =E\sum_{k=0}^{n\tau_n-1}\biggl[\Delta\biggl(\frac{k+1}{n}\biggr)
H(S_{k+1}-b_{k+1})-\Delta\biggl(\frac{k}{n}\biggr)
H(S_{k+1}-b_k)\biggr],
\end{eqnarray*}
with the last equality from \eqref{Heqn}, since $\{k<n\tau_n\}\in\F_k$.
The magnitude of the final expectation here is bounded by the sum of
%
\begin{equation}\label{excessterm1}
E\sum_{k=0}^{n\tau_n-1}\biggl|\Delta\biggl(\frac{k+1}{n}\biggr)
-\Delta\biggl(\frac{k}{n}\biggr)\biggr||H(S_{k+1}-b_{k+1})|
\end{equation}
and
%
\begin{equation}\label{excessterm2}
E\sum_{k=0}^{n\tau_n-1}\biggl|\Delta\biggl(\frac{k}{n}\biggr)\biggr|
|H(S_{k+1}-b_{k+1})-H(S_{k+1}-b_k)|.
\end{equation}
Using Corollaries~\ref{cor:sumbounds} and~\ref{cor:unifint}, it is
easy to show that~\eqref{excessterm1} tends to zero as $n\to\infty$.
To show that~\eqref{excessterm2} also tends to zero, we need a few
results from renewal theory and the fluctuation theory for random
walk.

Let $Y,Y_1,Y_2,\ldots$ be i.i.d. with $Y\sim S_{T_0}$, the first
ascending ladder height for $S_k$, $k\ge1$, and let $V_k=Y_1+\cdots+Y_k$
with $V_0\eqdef0$. Then $EY^3<\infty$ and the characteristic function
for $Y$ satisfies Cramer's condition. Define
\[
\mu(B)=\sum_{k=0}^\infty P(V_k\in B),
\]
so $\mu$ is the renewal measure for the random walk $V_k$, $k\ge0$.
By Wald's identity, for $x<0$,
\[
\mu((-\infty,-x))=ES_{T_x}/EY.
\]
So, for $\epsilon>0$ and $x<-\epsilon/2$,
%
\begin{equation}\label{Hdiff}
H(x+\epsilon/2)-H(x-\epsilon/2)
=\epsilon-EY\mu([-x-\epsilon/2,-x+\epsilon/2)).
\end{equation}
The following lemma follows immediately from these equations
and Theorem~3 of Stone~\cite{st65}.

\begin{lem}\label{lem:stone}
As $x\to-\infty$,
\[
H(x)=-\frac{E(Y+x)_+^2}{2(EY)^2}+\mathrm{o}\biggl(\frac{\log|x|}{x^2}\biggr).
\]
Also,
\[
H(x+\epsilon/2)-H(x-\epsilon/2)
=-\frac{\epsilon E(Y+x)_+}{EY}+\mathrm{o}\biggl(\frac{\log|x|}{|x|^3}\biggr),
\]
as $x\to-\infty$, uniformly for $\epsilon>0$ in any bounded set.
\end{lem}

To show that \eqref{excessterm2} is small we will need the following
lemma, similar to Lemma~\ref{lem:edge}, but bounding the expected
number of visits to smaller sets.

\begin{lem}\label{lem:smallsets}
Let $I_k=I_k(n)=(c_k,d_k)$ be intervals with
\[
\sup_{n,k}d_k\le K,
\]
and
\[
\sup_{n,k}(d_k-c_k)\le\frac{K}{\sqrt n},
\]
for some $K\in(0,\infty)$. If
\[
W=W_n=\#\{k\le n\tau_n\dvt b_k-S_k\in I_k\},
\]
then $EW\to0$ as $n\to\infty$.
\end{lem}

\begin{pf}
For $\alpha>0$, let
\[
W_0=\#\bigl\{k\le\min\{3\alpha n,n\tau_n\}\dvt b_k-S_k\in I_k\bigr\},
\]
the contribution to the count in $W$ from indices $k\le3\alpha n$.
Then $W_0\le M_{(-\infty,K]}(3\alpha)$. By Donsker's theorem,
\begin{eqnarray*}
\limsup_{n\to\infty}P\bigl(M_{(-\infty,K]}(3\alpha)\ge1\bigr)
&\le&\limsup_{n\to\infty} P\bigl(W_n(t)\ge b(t)-K/\sqrt n,\mbox{ for
some }t\le3\alpha\bigr)\\
&=&P\bigl(W(t)\ge b(t),\mbox{ for some }t\le3\alpha\bigr),
\end{eqnarray*}
which tends to $0$ as $\alpha\downarrow0$. So by Lemma~\ref{lem:edge},
$\limsup EW_0$ will be arbitrarily small if $\alpha$ is chosen suitably
small. Thus
this lemma will hold if
\[
EW-EW_0=\sum_{k>3\alpha n}P(b_k-S_k\in I_k,k\le n\tau_n)\to0,
\]
as $n\to\infty$ for any fixed $\alpha>0$.

Let $n^*=\lfloor\alpha n\rfloor$. By the local limit theorem of Stone~\cite{st65a}
(or Edgeworth expansion), for some~$K_0$,
%
\begin{equation}\label{local}
P(x-S_{n^*}\in I_k)\le K_0/n,
\end{equation}
for all $x\in\R$, $n\ge1$ and $k\ge1$. If $j<k$, $k-j<\sqrt n$,
$k\le n\tau_n$ and $b_k-S_k\in I_k$, then
\[
b_k-S_k\le K \quad\mbox{and}\quad  S_j<b_j.
\]
Together, these imply
\[
S_k-S_j\ge b_k-b_j-K\ge-\frac{k-j}{\sqrt n}\Vert b'\Vert_\infty-K
\ge-K-\Vert b'\Vert_\infty.
\]
Thus, if $n$ is large enough that $\alpha n>\sqrt n$ and if $k>3\alpha n$,
%
\begin{eqnarray}\label{smallsetbound}
&&
P(b_k-S_k\in I_k,k\le n\tau_n)
\nonumber
\\[-8pt]
\\[-8pt]
\nonumber
&&\quad\le P\bigl(b_k-S_k\in I_k ,n\tau_n\ge k/3,
S_k-S_j>-K-\Vert b'\Vert_\infty,0<k-j<\sqrt n\bigr).
\end{eqnarray}
To use this bound, let
\[
\G=\sigma\bigl(X_1,\ldots,X_{\lfloor k/3\rfloor},
X_{\lfloor k/3\rfloor+n^*+1},\ldots,X_k\bigr).
\]
Writing
\[
b_k-S_k=b_k-S_{\lfloor k/3\rfloor}-\bigl(S_k-S_{\lfloor k/3\rfloor
+n^*}\bigr)
-\bigl(S_{\lfloor k/3\rfloor+n^*}-S_{\lfloor k/3\rfloor}\bigr),
\]
since $b_k-S_{\lfloor k/3\rfloor}-(S_k-S_{\lfloor k/3\rfloor
+n^*})$
is $\G$ measurable and
\[
S_{\lfloor k/3\rfloor+n^*}-S_{\lfloor k/3\rfloor}\bigm|\G\sim S_{n^*},
\]
by \eqref{local},
\[
P(b_k-S_k\in I_k|\G)\le K_0/n.
\]
Since events $\{n\tau_n\ge k/3\}$ and
$\{S_k-S_j>-K-\Vert b'\Vert_\infty,0<k-j<\sqrt n\}$ are independent
and both lie in $\G$, using \eqref{smallsetbound} and conditioning on
$\G$,
\begin{eqnarray*}
 P(b_k-S_k\in I_k,k\le n\tau_n)
&\le&\frac{K_0}{n}P(n\tau_n\ge k/3)
P\bigl(S_k-S_j>-K-\Vert b'\Vert_\infty,
0<k-j<\sqrt n\bigr)\\
&=&\frac{K_0}{n}P(n\tau_n\ge k/3)
P\bigl(S_j>-K-\Vert b'\Vert_\infty,0<j<\sqrt n\bigr).
\end{eqnarray*}
The second probability in this bound tends to zero, and so
\[
EW-EW_0=\mathrm{o}(1)\frac1n\sum_{k>3\alpha n}P(n\tau_n>k/3)
=\mathrm{o}(1)E\tau_n\to0,
\]
proving the lemma.
\end{pf}

\begin{pf*}{Proof of Theorem~\protect\ref{thm:R}}
From the discussion and bounds above, the desired result will hold
if~\eqref{excessterm2} tends to zero, or if
\[
E\sum_{k=0}^{n\tau_n-1}|H(S_{k+1}-b_{k+1})-H(S_{k+1}-b_k)|\to0.
\]
The expectation of the final term in the sum tends to zero, for if
$S_{n\tau_n}-b_{n\tau_n-1}>0$, the summand is $|b_{n\tau_n}-b_{n\tau_n-1}|
\le\Vert b'\Vert_\infty/\sqrt n$; and if $S_{n\tau_n}-b_{n\tau_n-1}\le0$,
$S_{n\tau_n}-b_{n\tau_n}$ is within some multiple of $1/\sqrt n$ of
zero and the expectation will tend to zero by Lemma~\ref{lem:smallsets}.
So if $c_k=(b_{k}+b_{k-1})/2$ and
$\epsilon_k=|b_{k}-b_{k-1}|\le\Vert b'\Vert_\infty/\sqrt n$, we
need to show that
\[
E\sum_{k=1}^{n\tau_n-1}|H(S_k-c_k+\epsilon_k/2)
-H(S_k-c_k-\epsilon_k/2)|\to0.
\]
Using Lemma~\ref{lem:stone}, for some constant $K_0$,
\[
|H(x+\epsilon/2)-H(x-\epsilon/2)|\le\epsilon g(x)
+\frac{K_0}{|x|^{5/2}},
\]
for all $x<0$ and all $\epsilon\in[0,\Vert b'\Vert_\infty]$, where
$g(x)=E(Y+x)_+/EY$.
Using this,
\begin{eqnarray*}
&& E\sum_{k=1}^{n\tau_n-1}|H(S_k-c_k+\epsilon_k/2)
-H(S_k-c_k-\epsilon_k/2)|I\{b_k-S_k\ge K_1\}\\
&&\quad\le\frac{\Vert b'\Vert_\infty}{\sqrt n}E\sum_{k=1}^{n\tau_n-1}
g(S_k-c_k)
+K_0E\sum_{k=1}^{n\tau_n-1}\min\{K_1^{-5/2},(b_k-S_k)^{-5/2}\}.
\end{eqnarray*}
By Corollary~\ref{cor:sumbounds},
\begin{eqnarray*}
E\sum_{k=1}^{n\tau_n-1}g(S_k-c_k)
&\le& K\biggl[\frac{E(Y+\Lambda)_+}{EY}+\frac{2}{EY}
E\int|x|(Y+x+\Lambda)_+I\{x<0\} \,\mathrm{d}x\biggr]\\
&=&K\biggl[\frac{E(Y+\Lambda)_+}{EY}+\frac{E(Y+\Lambda)^3}{3EY}\biggr],
\end{eqnarray*}
which is finite since $EY^3<\infty$. And by the same corollary,
\begin{eqnarray*}
E\sum_{k=1}^{n\tau_n-1}
\min\{K_1^{-5/2},(b_k-S_k)^{-5/2}\}
&\le& K\biggl[K_1^{-5/2}+2\int^0|x|\min\{K_1^{-5/2},|x|^{-5/2}\}
 \,\mathrm{d}x\biggr]\\
&=&K\biggl[K_1^{-1/2}+\frac43K_1^{-3/2}+K_1^{-5/2}\biggr].
\end{eqnarray*}
Since this bound tends to zero as $K_1\to\infty$, the theorem will
hold if
\[
E\sum_{k=1}^{n\tau_n-1}|H(S_k-c_k+\epsilon_k/2)
-H(S_k-c_k-\epsilon_k/2)|I\{b_k-S_k<K_1\}
\]
converges to zero for any fixed $K_1$. Also, using Lemma~\ref{lem:smallsets},
we can include the restriction $b_k-S_k>\Vert b'\Vert_\infty/\sqrt n$ in
the indicator. Using \eqref{Hdiff}, it will then be sufficient to show
\[
E\sum_{k=1}^{n\tau_n-1}\epsilon_kI\{b_k-S_k<K_1\}\to0
\]
and
\[
E\sum_{k=1}^{n\tau_n-1}\mu([c_k-S_k-\epsilon_k/2,
c_k-S_k+\epsilon_k/2])I\biggl\{\frac{\Vert b'\Vert_\infty}{\sqrt n}
<b_k-S_k<K_1\biggr\}\to0.
\]
The first of these follows immediately from Lemma~\ref{lem:edge}.
Using Fubini's theorem, the second expression equals
\[
\int E\sum_{k=1}^{n\tau_n-1}I\biggl\{\frac{\Vert b'\Vert_\infty}{\sqrt n}
<b_k-S_k<K_1,|S_k-c_k+x|\le\epsilon_k/2\biggr\} \,\mathrm{d}\mu(x).
\]
By Lemma~\ref{lem:smallsets}, the integrand here tends to zero, uniformly
in $x$. Since the range of integration remains bounded, the integral must
tend to zero, proving the theorem.
\end{pf*}

\section*{Acknowledgements}
I would like to thank the referee for a
careful reading of this paper and several useful suggestions.
Supported in part by NSA Grant F012499 and NSF Grant 0706771.


\printhistory

\end{document}